\documentclass[a4paper]{article}
\usepackage[margin=30mm]{geometry}
\usepackage{amsmath}
\usepackage{amssymb}
\usepackage{amsthm}
\usepackage{braket}
\usepackage{titlesec}
\usepackage{titlefoot}
\usepackage{comment}
\usepackage{bm}
\titleformat*{\section}{\normalsize\bfseries}
\titleformat*{\subsection}{\normalsize\bfseries}
\allowdisplaybreaks

\newtheorem{thm}{Theorem}

\newtheorem{rem}[thm]{Remark}

 
\begin{document}

\title{
\vspace{-1.25cm}
\Large{\bf Higher-order Asymptotic Expansions for \\
Laplace's Integral and Their Error Estimates}}
\author{Ikki Fukuda and Yoshiki Kagaya
}
\date{}
\maketitle

\footnote[0]{This is the preprint version of the original paper, including its erratum, published in ``The 43rd JSST Annual International Conference on Simulation Technology, Conference Proceedings (2024) 502--508''.}

\vspace{-1cm}
\begin{abstract}
We deal with the asymptotic analysis for Laplace's integral. 
For this problem, the so-called ``Laplace's method'' by P.S.~Laplace (1812) is well-known and it has been developed in various forms over many years of studies.
In this paper, we derive some formulas of the higher-order asymptotic expansions for that integral, with error estimates, which generalize previous results. 
Moreover, we discuss a comparison of these asymptotic formulas with approximations using numerical integral.
\end{abstract}

\medskip
\noindent
{\bf Keywords:} 
Laplace's integral; Higher-order asymptotic expansions. 

\section{Introduction}  

We consider the asymptotic behavior of the following integral as $t\to \infty$:  
\begin{equation}
I(t)=\int_{a}^{b} e^{th(x)}g(x)dx, \ \ t>0,  \label{LI}
\end{equation}
where $a, b\in \mathbb{R}$ and the functions $g,h: [a,b] \to \mathbb{R}$ are continuous. 
This integral \eqref{LI} is usually called ``Laplace's integral'', which has several applications in various fields of science and engineering as well as mathematics, and the literature on this integral is quite extensive. For example, Laplace's integral \eqref{LI} often appears in the probability theory and statistics (see, e.g.~\cite{B22,B86,DZ98,K20,M82,DZ03,W19,WU24} and also references therein). 
In particular, we note that the asymptotic analysis for \eqref{LI} can be applied to the large deviation principle (cf.~\cite{DZ98,DZ03}). 
Moreover, it is applied to the study of particle filters that appear in the fields of signal processing and statistics \cite{MBL11}. 
In addition, it also has applications in the field of random chaos \cite{KPH15}. 
Therefore, we believe that our study of \eqref{LI} will be useful in not only pure mathematics, but also in many areas of science and engineering.

For the asymptotic behavior of the above integral \eqref{LI} as $t\to \infty$, the so-called ``Laplace's method'' or ``Laplace's approximation'' by Laplace~\cite{L12} is well-known and it has been developed in various forms over many years of studies. For example, we can refer to \cite{K20,N20,S00,W19,WU24} and also references therein. 
Here, we shall explain about the results of that asymptotic analysis, based on the textbook by Shiga~\cite{S00}.
Assume that $h(x)$ has a maximum only at $x=c \in [a, b]$, and is twice continuously differentiable function around $x=c$ satisfying $h''(c) \neq 0$. Then, if $g(c)\neq0$, $I(t)$ satisfies the following asymptotic formulas:

\medskip
\noindent
{\rm $( \mathrm{i} )$} If $a<c<b$, then 
\begin{equation}
I(t) \sim \sqrt{\frac{2\pi}{\left\lvert h''(c)\right\rvert } }
g(c)\frac{e^{th(c)}}{\sqrt{t}} \ \ (t \to \infty). \label{thm-Laplace-1}
\end{equation}
{\rm $( \mathrm{ii} )$} If $c=a$ or $c=b$, and $h'(c)=0$, then
\begin{equation}
I(t) \sim \sqrt{\frac{\pi}{2\left\lvert h''(c)\right\rvert } }
g(c)\frac{e^{th(c)}}{\sqrt{t}} \ \ (t \to \infty). \label{thm-Laplace-2}
\end{equation}
{\rm $( \mathrm{iii} )$} If $c=a$ or $c=b$, and $h'(c)\neq 0$, then 
\begin{equation}
I(t) \sim \frac{g(c)}{\left\lvert h'(c)\right\rvert} 
\frac{e^{th(c)}}{\sqrt{t}} \ \ (t \to \infty). \label{thm-Laplace-3}
\end{equation}
The notation ``$\sim$'' is used to mean that the quotient of the left hand side by the right hand side converges 1 as $t\to \infty$. 
As we mentioned in the above, some related results to these formulas have already been obtained by many researchers. 
Especially, as a result of its most generalized form up to the authors knowledge, let us refer to the recent paper by Nemes~\cite{N20}. 

\newpage

Although the above results \eqref{thm-Laplace-1}, \eqref{thm-Laplace-2} and \eqref{thm-Laplace-3} are very useful, there are some problems from the perspective of application. 
For example, when the maximum point of $h(x)$ and the zero point of $g(x)$ overlap, i.e. when $g(c)=0$, these formulas do not make sense, 
because the asymptotic profiles vanish identically. 
Therefore, for the more general $I(t)$ in which such situations occur, we need to construct some new asymptotic formulas. 
Moreover, from application point of view, the errors between $I(t)$ and the asymptotic profiles should also be investigated.
In our study, we analyzed Laplace's integral \eqref{LI} based on this consideration and succeeded to derive some results on the higher-order asymptotic expansions for $I(t)$, with error estimates, which overcomes this situation.

\section{Main Result}

Now, let us state our main result in this paper: 
\begin{thm}\label{main-ap}
Assume that $h(x)$ has a maximum only at $x=c \in [a, b]$, and is three times continuously differentiable function on $[a, b]$ satisfying $h''(c) \neq 0$. 
Moreover, let $k\in \mathbb{N}\cup \{0\}$ and suppose that $g(x)$ is $(k+1)$-times continuously differentiable function on $[a, b]$ satisfying $g(c) = g'(c) = \cdots =g^{(k-1)}(c) = 0$ and $g^{(k)}(c) \neq 0$. Then, $I(t)$ satisfies the following asymptotic formulas:

\medskip
\noindent
{\rm $( \mathrm{I} )$} If $a < c < b$, we additionally assume that $k$ is even number. Then, 
\begin{equation}
I(t) = e^{th(c)} \left( t^{-\frac{1}{2}-\frac{k}{2}}
g^{(k)}(c) \sqrt{\frac{2\pi}{{\left\lvert h''(c) \right\rvert}^{k+1}}}\frac{(k-1)!!}{k!} + O(t^{-1-\frac{k}{2}})\right)
\ \ (t \to \infty). \label{main-1}
\end{equation}
{\rm $\rm(I\hspace{-.15em}I)$} If $c=a$ or $c=b$, and $h'(c)=0$, then
\begin{equation}
I(t) = e^{th(c)} \left( t^{-\frac{1}{2}-\frac{k}{2}}
g^{(k)}(c) \sqrt{\frac{2^{k-1}}{\left\lvert h''(c) \right\rvert^{k+1} }}\,\gamma(k, c)+ O(t^{-1-\frac{k}{2}})\right)
\ \ (t \to \infty), \label{main-2}
\end{equation}
where the constant $\gamma(k, c)$ is defined by 
\[
\gamma(k, c):=
\begin{cases}
\displaystyle \sqrt{\frac{\pi}{2^{k}}} \frac{(k-1)!!}{k!} & (k:\text{even}), \\[.8em]
\displaystyle \frac{1}{k!} \left( \frac{k-1}{2} \right)! & (k:\text{odd}, \ c=a), \\[.8em]
\displaystyle -\frac{1}{k!} \left( \frac{k-1}{2} \right)! & (k:\text{odd}, \ c=b). 
\end{cases}
\]
{\rm $\rm(I\hspace{-.15em}I\hspace{-.15em}I)$} If $c=a$ or $c=b$, and $h'(c)\neq0$, then
\begin{equation}
I(t) =
e^{th(c)} \left(t^{-1-k} g^{(k)}(c) (h'(c))^{-1-k} (-1)^{k}\eta(c)+O(t^{-2-k})\right) \ \ (t \to \infty), \label{main-3}
\end{equation}
where the constant $\eta(c)$ is defined by $\eta(c):=-1$ if $c=a$ and $\eta(c):=1$ if $c=b$.
\end{thm}

\begin{rem}
For some related results on the asymptotic formulas of $I(t)$, especially with error estimates, have already been obtained in e.g. Wakaki~{\rm \cite{W19}} and Wakaki--Ulyanov~{\rm \cite{WU24}} (see, also Kolokoltsov~{\rm \cite{K20}}). 
Our formulas can be considered as one of the generalization of them including the higher-order asymptotic expansion.
\end{rem}

\section{Proof of the Main Result}

In this section, we shall prove our main result Theorem~\ref{main-ap}. 
Here, we would like to give its proof by applying the techniques used in \cite{W19,WU24}. 
In particular, we note that it is given by the direct calculations based on Taylor's expansion. 

\medskip
\noindent{\it {\bf Proof of Theorem~\ref{main-ap}.}} 
We shall only prove the case of {\rm $( \mathrm{I} )$}, since the other two cases can be proved by slightly modifying the technique of {\rm $( \mathrm{I} )$}. 
First, from the assumptions on $h(x)$ and $a<c<b$, we note that $h(x)$ satisfies $h'(c) = 0$ and $h''(c) < 0$. 
Moreover, since $h(x)$ is three times continuously differentiable on $[a,b]$, applying Taylor's theorem to $h(x)$, there exists $\theta_{0} \in (0,1)$ such that 
\newpage
\begin{equation}
h(x) = h(c) + \frac{1}{2} h''(c)(x-c)^2 + \frac{1}{3!} h^{(3)}(c + \theta_{0}(x-c))(x-c)^3.  \label{L1_Taylor}
\end{equation}

In what follows, let us consider the asymptotic expansion for $I(t)e^{-th(c)}$. 
First, we take sufficiently small $\delta>0$. Now, we split the integral and define the functions $J_{1}(t)$ and $J_{2}(t)$ as follows: 
\begin{align}
I(t) e^{-th(c)} 
&=\left(\int_{c-\delta}^{c+\delta} +\int_{[a,b] \backslash [c-\delta,c+\delta]} \right)e^{t(h(x) - h(c))} g(x) dx \nonumber \\
&=:J_{1}(t)+J_{2}(t). \label{Jt_def}
\end{align}
First, we would like to rewrite $J_1(t)$. Now, noticing the assumptions on $g(x)$ and applying Taylor's theorem to $g(x)$, then there exists $\theta_1 \in (0,1)$ such that 
\begin{equation}
g(x) = \frac{g^{(k)}(c)}{k!} (x-c)^{k}+ \frac{1}{(k+1)!} g^{(k+1)}(c + \theta_1(x-c))(x-c)^{k+1}. \label{g-Taylor}
\end{equation}
We set $\sigma := \sqrt{\left\lvert h''(c) \right\rvert} $ for simplicity. Then, it follows from \eqref{L1_Taylor} and \eqref{g-Taylor} that 
\begin{align*}
	J_1(t) 
	&= \int_{c-\delta}^{c+\delta} \exp \left\{ t \left( \frac{1}{2} h''(c) (x-c)^2 
	+ \frac{1}{3!} h^{(3)}(c+\theta_{0}(x-c))(x-c)^3 \right) \right\} \\
	&\ \ \ \ \times \left\{ \frac{g^{(k)}(c)}{k!} (x-c)^{k}
	+ \frac{1}{(k+1)!} g^{(k+1)}\left(c + \theta_{1}(x-c)\right) (x-c)^{k+1} \right\} dx \\
	&=\sigma^{-1}t^{-\frac{1}{2}}\int_{-\delta \sigma \sqrt{t}}^{\delta \sigma \sqrt{t}} e^{-\frac{y^2}{2}}
\exp \left\{\frac{1}{3!}h^{(3)}\left(c+\theta_{0}\frac{y}{\sigma \sqrt{t}}\right)\frac{1}{\sqrt{t}}\left(\frac{y}{\sigma}\right)^3\right\} \nonumber \\
&\ \ \ \ \times \left\{\frac{g^{(k)}(c)}{k!}\left(\frac{y}{\sigma\sqrt{t}}\right)^{k}
+ \frac{1}{(k+1)!} g^{(k+1)}\left(c + \theta_1\frac{y}{\sigma\sqrt{t}}\right)\left(\frac{y}{\sigma\sqrt{t}}\right)^{k+1}\right\} dy.  \nonumber
	\end{align*}
Moreover, applying Taylor's theorem to the exponential function in the above integral, there exists $\theta_2 \in (0,1)$ such that the following relation holds: 
\begin{align}
J_{1}(t)
&= \frac{g^{(k)}(c)}{k!} \sigma^{-1-k}t^{-\frac{1}{2}-\frac{k}{2}}\int_{-\delta \sigma \sqrt{t}}^{\delta \sigma \sqrt{t}} 
e^{-\frac{y^2}{2}}y^{k}dy \nonumber \\
&\ \ \ \ +\frac{ \sigma^{-2-k}}{(k+1)!}t^{-1-\frac{k}{2}}\int_{-\delta \sigma \sqrt{t}}^{\delta \sigma \sqrt{t}}
g^{(k+1)}\left(c + \theta_1\frac{y}{\sigma\sqrt{t}}\right) e^{-\frac{y^2}{2}} y^{k+1}dy \nonumber \\ 
&\ \ \ \ +
\frac{\sigma^{-4}}{3!}t^{-1}\int_{-\delta \sigma \sqrt{t}}^{\delta \sigma \sqrt{t}}
h^{(3)}\left(c+\theta_{0}\frac{y}{\sigma \sqrt{t}}\right)y^3 \nonumber \\
&\quad \quad \ \times \exp \left\{ -\frac{y^2}{2} + \frac{\theta_{2}}{3!}h^{(3)}\left(c+\theta_{0}\frac{y}{\sigma \sqrt{t}}\right)\frac{1}{\sqrt{t}}\left(\frac{y}{\sigma}\right)^3\right\} \nonumber \\
& \quad \quad \ \times \left\{\frac{g^{(k)}(c)}{k!}\left(\frac{y}{\sigma\sqrt{t}}\right)^{k}
+ \frac{1}{(k+1)!} g^{(k+1)}\left(c + \theta_1\frac{y}{\sigma\sqrt{t}}\right)\left(\frac{y}{\sigma\sqrt{t}}\right)^{k+1}\right\} dy \nonumber \\
&=: J_{1.1}(t) + J_{1.2}(t)+J_{1.3}(t). \label{main_1.1_J1y}
\end{align}

Next, we shall derive the leading term of $I(t)e^{-th(c)}$. In order to do that, let us transform the above $J_{1.1}(t)$ as follows: 
\begin{align}
J_{1.1}(t)
&= \frac{g^{(k)}(c)}{k!} \sigma^{-1-k}t^{-\frac{1}{2}-\frac{k}{2}}\left(\int_{-\infty}^{\infty}-\int_{|y| \geq \delta \sigma \sqrt{t}}\right) e^{-\frac{y^2}{2}} y^{k}dy  \nonumber \\
&=: J_{1.1.1}(t) - J_{1.1.2}(t). \label{split-J1.1}
\end{align}
Then, the leading term can be found from $J_{1.1.1}(t)$. Actually, since $k$ is even number, it follows from a direct calculation that 
\begin{equation}
J_{1.1.1}(t)
=t^{-\frac{1}{2}-\frac{k}{2}} g^{(k)}(c) \sqrt{\frac{2\pi}{{\left\lvert h''(c) \right\rvert}^{k+1}}}\frac{(k-1)!!}{k!}. \label{L1c_M}
\end{equation}
On the other hand, we are able to consider $J_{1.1.2}(t)$ as a remainder term. Indeed, we can easily see that the following exponential decay estimate holds: 
\begin{align}
&\left|J_{1.1.2}(t)\right|\le Ct^{-\frac{1}{2}-\frac{k}{2}}  \left\lvert \int_{|y| \geq \delta \sigma \sqrt{t}}e^{-\frac{y^2}{2}} y^{k}dy \right\rvert \nonumber \\
& \leq Ct^{-\frac{1}{2}-\frac{k}{2}} \sup_{|y| \geq \delta \sigma \sqrt{t}} e^{-\frac{y^2}{4}}\int_{-\infty}^{\infty}e^{-\frac{y^2}{4}} \left|y\right|^{k}dy 
= O\left(t^{-\frac{1}{2}-\frac{k}{2}} e^{-\frac{\delta^2\sigma^2t}{4}}\right). \label{L1c_J1.1.1}
\end{align}

In the rest of this proof, let us evaluate the other two remainder terms $J_{1.2}(t)$ and $J_{1.3}(t)$. 
First, in the similar way to get \eqref{L1c_J1.1.1}, we obtain 
\begin{align}
\left\lvert J_{1.2}(t) \right\rvert
& \leq Ct^{-1-\frac{k}{2}} \displaystyle \sup_{x \in [a,b]}\left|g^{(k+1)}(x)\right| 
\int_{-\infty}^{\infty} e^{-\frac{y^2}{2}}\left|y\right|^{k+1} dy =
O(t^{-1-\frac{k}{2}}). \label{L1c_J1.1.2}
\end{align}
Next, we shall treat $J_{1.3}(t)$. Before evaluating it, let us rewrite it as follows: 
\begin{align}
J_{1.3}(t)
&=\frac{g^{(k)}(c)}{3!k!}\sigma^{-4-k}t^{-1-\frac{k}{2}}\int_{-\delta \sigma \sqrt{t}}^{\delta \sigma \sqrt{t}}
h^{(3)}\left(c+\theta_{0}\frac{y}{\sigma \sqrt{t}}\right)y^{3+k} \nonumber \\
& \ \ \ \ \times \exp \left\{ -\frac{y^2}{2} + \frac{\theta_{2}}{3!} h^{(3)}\left(c+\theta_{0}\frac{y}{\sigma \sqrt{t}}\right)\frac{1}{\sqrt{t}}\left(\frac{y}{\sigma}\right)^3\right\} dy \nonumber \\
& \ \ \ \ +\frac{\sigma^{-5-k}}{3!(k+1)!}t^{-\frac{3}{2}-\frac{k}{2}}\int_{-\delta \sigma \sqrt{t}}^{\delta \sigma \sqrt{t}}
h^{(3)}\left(c+\theta_{0}\frac{y}{\sigma \sqrt{t}}\right) g^{(k+1)}\left(c+\theta_{1}\frac{y}{\sigma \sqrt{t}}\right)y^{4+k} \nonumber \\
& \ \ \ \ \times \exp \left\{ -\frac{y^2}{2} + \frac{\theta_{2}}{3!} h^{(3)}\left(c+\theta_{0}\frac{y}{\sigma \sqrt{t}}\right)\frac{1}{\sqrt{t}}\left(\frac{y}{\sigma}\right)^3\right\} dy.
\end{align}
Now, we would like to evaluate the exponential function in the above. 
Since, $-\delta \sigma \sqrt{t} \leq y \leq \delta \sigma \sqrt{t} \Leftrightarrow \lvert \displaystyle y/(\sigma \sqrt{t}) \rvert  \leq \delta$, we have 
\begin{align*}
&\exp \left\{ -\frac{y^2}{2} + \frac{\theta_{2}}{3!}h^{(3)}\left(c+\theta_{0}\frac{y}{\sigma \sqrt{t}}\right)\frac{1}{\sqrt{t}}\left(\frac{y}{\sigma}\right)^3\right\} \\
&\leq \exp \left\{ -\frac{y^2}{2} \left( 1 - \displaystyle \sup_{x \in [a,b]}\left|h^{(3)}(x)\right| \frac{\delta}{\sigma^2} \right) \right\}
=: e^{-\gamma y^2}. 
\end{align*}
Noticing that $\gamma>0$ if we choose $\delta>0$ sufficiently small. Thus, we can see that 
\begin{align}
\left\lvert J_{1.3}(t) \right\rvert
& \leq Ct^{-1-\frac{k}{2}} \sup_{x \in [a,b]}\left|h^{(3)}(x)\right| 
\int_{-\infty}^{\infty} e^{-\gamma y^2}\left|y\right|^{3+k} dy   \nonumber \\
&\ \ \ +Ct^{-\frac{3}{2}-\frac{k}{2}}  \sup_{x \in [a,b]}\left|h^{(3)}(x)\right| \displaystyle \sup_{x \in [a,b]}\left|g^{(k+1)}(x)\right|\int_{-\infty}^{\infty} e^{-\gamma y^2}\left|y\right|^{4+k} dy \nonumber \\
& = O(t^{-1-\frac{k}{2}}). \label{L1c_J1.2.2}
\end{align}

Finally, we would like to deal with $J_{2}(t)$. Now, we set 
\begin{equation*}
\alpha := - \sup \left\{ h(x)-h(c) : x\in [a,b] \backslash [c-\delta, c+\delta] \right\}>0. 
\end{equation*}
Then, we can easily evaluate $J_{2}(t)$ in \eqref{Jt_def} as follows: 
\begin{align}
\left\lvert J_2(t) \right\rvert 
\leq e^{-\alpha t} \int_{[a,b] \backslash [c-\delta,c+\delta]} \left\lvert g(x) \right\rvert dx
\leq e^{-\alpha t} \int_{a}^{b} \left\lvert g(x) \right\rvert dx = O(e^{-\alpha t}). \label{L1c_J2}
\end{align}
Eventually, combining \eqref{Jt_def}, \eqref{main_1.1_J1y} through \eqref{L1c_J1.1.2}, \eqref{L1c_J1.2.2} and \eqref{L1c_J2}, we can conclude that the desired result \eqref{main-1} is true. 
This completes the proof of {\rm $( \mathrm{I} )$}. \qed

\section{Comparison with a Numerical Integral}

When investigating the asymptotic behavior of $I(t)$ as $t\to \infty$, one may be more familiar with simulations based on numerical integrals rather than deriving rigorous asymptotic profiles.
In the rest of this paper, we discuss the relationship between our asymptotic formulas \eqref{main-1}, \eqref{main-2} and \eqref{main-3} for $I(t)$, and approximation of $I(t)$ by using a famous numerical integral. 
First, let us recall some known results about that numerical integral. Here, we consider the following integral: 
\begin{equation*}
I=\int_{a}^{b}f(x)dx. 
\end{equation*}
A well-known method for approximating this integral is a numerical integral called the ``composite Simpson's rule''. 
Now, we shall explain about this method briefly. 
Let $n\ge2$ be an even number. We split the interval $[a, b]$ into $n$-th parts with the mesh size $h$ (i.e. $h:=(b-a)/n$) and set $x_{i}:=a+ih$ ($i=0, 1, \cdots, n$). 
Then, the above integral can be approximated by 
\begin{equation*}
S_{n}[f]:=\frac{h}{3}\left(f(x_{0})+2\sum_{i=1}^{n/2-1}f(x_{2i})+4\sum_{i=1}^{n/2}f(x_{2i-1})+f(x_{n})\right). 
\end{equation*}
It is known that the maximum error bound of the approximation by using this method can be evaluated as follows (for details, see e.g.~\cite{G12}):
\begin{equation*}
\left|I-S_{n}[f]\right|\le \frac{(b-a)^{5}}{180n^{4}}\max_{\xi \in [a, b]}\left|f^{(4)}(\xi)\right|. 
\end{equation*}

Here, applying the composite Simpson's rule to Laplace's integral $I(t)$ in \eqref{LI}, as $f_{t}(x)=e^{th(x)}g(x)$ with the parameter $t>0$, then it follows from the above error estimate that 
\begin{equation}\label{simpson-error}
\left|I(t)-S_{n}[f_{t}]\right|\le \frac{c_{0}}{n^{4}}\exp\left(t\max_{\xi \in [a, b]}\left|h(\xi)\right|\right)t^{4}=:E_{S}(t), \ \ t\gg1, 
\end{equation}
where $c_{0}>0$ is a certain positive constant. On the other hand, for example in the case of {\rm (I)} in Theorem~\ref{main-ap}, according to our result \eqref{main-1}, there exists a positive constant $C_{0}>0$ such that 
\begin{align}
&\left|I(t)-e^{th(c)} t^{-\frac{1}{2}-\frac{k}{2}}
g^{(k)}(c) \sqrt{\frac{2\pi}{{\left\lvert h''(c) \right\rvert}^{k+1}}}\frac{(k-1)!!}{k!}\right| \nonumber \\
&\le C_{0}\exp\left(t\max_{\xi \in [a, b]}\left|h(\xi)\right|\right)t^{-1-\frac{k}{2}}=:E_{L}(t), \ \ t\gg1. \label{our-result}
\end{align}
Therefore, compared with \eqref{simpson-error} and \eqref{our-result}, we can see that  
\[
E_{L}(t)<E_{S}(t) \ \Leftrightarrow \ C_{0}t^{-1-\frac{k}{2}}<\frac{c_{0}}{n^{4}}t^{4} \ \Leftrightarrow \ t>\left(\frac{C_{0}n^{4}}{c_{0}}\right)^{\frac{1}{5+\frac{k}{2}}}. 
\]
Thus, we can conclude that for any mesh number $n$, there exists a sufficiently large $T>0$ such that $E_{L}(t)<E_{S}(t)$ holds if $t>T$. 
In other words, our results allow us to approximate $I(t)$ so that the maximum error is smaller than the maximum error in the composite Simpson's rule, for $t\gg1$. 
However, we note that the asymptotic formulas do not necessarily give a good approximation when $t>0$ is small. 
Therefore, we think that both numerical calculations and asymptotic formulas should be used appropriately depending on the situation.

\section*{Acknowledgments}

This study is supported by Grant-in-Aid for Young Scientists Research No.22K13939, Japan Society for the Promotion of Science. The authors would like to express their appreciation to Professor Yuki Ueda in Hokkaido University of Education for his useful comments and stimulating discussions.

\section*{Disclosure of Interests}

The authors have no competing interests to declare that are relevant to the content of this article.



\medskip
\par\noindent
\begin{flushleft}
Ikki Fukuda\\
Faculty of Engineering, Shinshu University, \\
Nagano, 380-8553, Japan\\
E-mail: i\_fukuda@shinshu-u.ac.jp

\bigskip
\par\noindent
Yoshiki Kagaya\\
Department of Mathematical Science, Shinshu University, \\
Matsumoto, 390-8621, Japan\\
23ss104g@shinshu-u.ac.jp
\end{flushleft}

\newpage

The following sentences are part of the erratum issued for the paper: 

\section*{Addendum}

After the publication of the original paper, it was pointed out to us that our main result Theorem 1 can be regarded as direct or almost immediate consequences of Theorem 2.1 in Chapter 9, Section 2 of \cite{O74} by F.W.J. Olver. Actually, the case {\rm $\rm(I\hspace{-.15em}I)$} in Theorem 1 deals with Laplace's integral where the phase function $h(x)$ attains a maximum at an endpoint of the integration interval, and the integrand behaves like $(x-c)^{k}$ near that point. This situation matches the setting treated in \cite{O74}, which provides asymptotic expansions for such integrals under more general conditions. Moreover, the case {\rm $\rm(I)$} can be derived by splitting the integral at the critical point and applying the formula derived in {\rm $\rm(I\hspace{-.15em}I)$} to each resulting piece. Furthermore, the case {\rm $\rm(I\hspace{-.15em}I\hspace{-.15em}I)$} corresponds to the special case of Olver's result with $\mu=1$ in Theorem 2.1 of \cite{O74}. Although these connections were not explicitly mentioned in the original manuscript, we hope this clarification helps position the results within the broader literature. We believe that the paper remains a useful reference, particularly due to its explicit derivations, concrete error estimates, and comparisons with numerical integration, which may assist readers applying this method in practical science and engineering. Finally, the authors would like to express their sincere gratitude to Professor Gerg\H{o} Nemes for bringing this to our attention and providing a lot of valuable comments. 


\begin{thebibliography}{99}
\addcontentsline{toc}{section}{References}

\bibitem{B22}
Bardet, J.-M.: Laplace's method and BIC model selection for least absolute value criterion. Statistics and Probability Letters {\bf 195}, 109764, 10 pp. (2023).

\bibitem{B86}
Bolthausen, E.: Laplace approximations for sums of independent random vectors. Probab. Th. Rel. Fields {\bf 72}, 305--318 (1986). 

\bibitem{DZ98}
Dembo, A., Zeitouni, O.: Large deviations techniques and applications (Second Edition). Springer-Verlag, New York (2009).

\bibitem{G12}
Gautschi, W.: Numerical Analysis (Second Edition). Birkh\"{a}user, Boston (2012). 

\bibitem{K20}
Kolokoltsov, V.N.: Rates of convergence in Laplace's integrals and sums and conditional central limit theorems. Mathematics {\bf 8}, 479, 19 pp. (2020).

\bibitem{KPH15}
Korshunov, D.A., Piterbarg, V.I., Hashorva, E.: On the asymptotic Laplace method and its application to random chaos. Mathematical Notes {\bf 97}, 878--891 (2015).

\bibitem{L12}
Laplace, P.S.: Th\'{e}orie analytique des probabilit\'{e}s. Ve. Courcier, Paris (1812). 

\bibitem{M82}
Martin-L\"{o}f, A.: Laplace approximations for sums of independent random variables. Z. Wahrscheinlichkeitstheorie verw. Gebiete {\bf 59}, 101--115 (1982). 

\bibitem{DZ03}
Moral, P.D., Zajic, T.: A note on the Laplace--Varadhan integral lemma. Bernoulli {\bf 9}, 49--65 (2003).

\bibitem{MBL11}
Musso, C., Quang, P.B., Gland, F.L.: Introducing the Laplace approximation in particle filtering. In: 14th International Conference on Information Fusion, pp. 290--297, Chicago, Illinois (2011).

\bibitem{N20}
Nemes, G.: An extension of Laplace's method. Constr. Approx. {\bf 51}, 247--272 (2020).

\bibitem{O74}
Olver, F.W.J.: Asymptotics and Special Functions. Academic Press, New York (1974).

\bibitem{S00}
Shiga, T.: Probability theory from Lebesgue integral, Kyoritsu-Shuppan, Tokyo (2000) (in Japanese).

\bibitem{W19}
Wakaki, H.: Laplace approximation and its applications. RIMS K${\rm \hat{o}}$ky${\rm \hat{u}}$roku {\bf 2133}, 66--74 (2019) (in Japanese).

\bibitem{WU24}
Wakaki, H., Ulyanov, V.V.: Laplace expansion for Bartlett--Nanda--Pillai test statistic and its error bound. Theory Prob. Appl. {\bf 68}, 570--581 (2024). 

\end{thebibliography}
\end{document}